\documentclass[12pt,thmsa]{article}
\usepackage{amsfonts}
\usepackage{amssymb}
\usepackage{amsmath,amsthm}
\usepackage[]{latexsym}
\newtheorem{teo}{Theorem}[section]

\newtheorem{obs2}[teo]{Remark}

\newtheorem{tea}{Theorem}[subsection]

\newtheorem{no2}[teo]{Note}

\newtheorem{no3}[tea]{Note}

\newcommand{\Gal}{{\rm Gal}}
\newcommand{\Fr}{{\rm Frob }}
\newcommand{\tr}{{\rm trace}}
\newcommand{\Q}{\mathbb{Q}}

\newcommand{\la}{\lambda}

\newcommand{\re}{\rho_\lambda}
\newcommand{\rer}{\bar{\rho}_\lambda}
\newcommand{\Z}{\mathbb{Z}}
\newcommand{\tres}{t}

\newcommand{\disc}{{\rm disc}}

\newcommand{\PGL}{{\rm PGL}}

\newcommand{\GL}{{\rm GL}}

\newcommand{\End}{{\rm End}}

\newcommand{\R}{{\mathbb R}}
\newcommand{\F}{{\mathbb F}}

\newcommand{\cO}{{\mathcal O}}

\title{Modularity of Abelian Surfaces with Quaternionic Multiplication}
\author{Luis  Dieulefait\\
Institut Galil{\'e}e, Universit{\'e}  Paris 13\\
99, av. J.B. Cl{\'e}ment\\
93430 Villetaneuse, France\\
e-mail: luisd@math.univ-paris13.fr, luisd@mat.ub.es}

\begin{document}

\maketitle

\vspace{6mm}

\newpage

\begin{abstract} We prove that any abelian surface defined over $\Q$ of
$GL_2$-type having quaternionic multiplication and good reduction at
$3$ is modular. We  generalize the result to
higher dimensional abelian varieties with ``sufficiently many
endomorphisms"
%Mathematics Subject Classification:
\end{abstract}

\section{Statement of the theorem}

In this brief article we interest ourselves in the relation to classical
modular forms of the following geometric object: $A$ an abelian
surface defined over $\Q$ such that $\End^0_\Q(A)= \Q(\sqrt{d})$ (d a non-square
integer)
 and
$\End^0_{\overline{\Q}}(A)= B$ where $B/\Q$ is an indefinite
rational quaternion algebra. \\
%For simplicity we will also assume that $\End_{\overline{\Q}}(A)$
 % is the maximal order in $B$.\\
 We borrow from P. Clark's note (see [C])
  the name ``premodular"  QM-surfaces over $\Q$ for
  these abelian surfaces. The condition of having real multiplication
  defined over $\Q$ is necessary (and sufficient)
  in order to obtain a two dimensional
  Galois representation of the full $\Gal(\overline{\Q}/\Q)$ on the
  Tate modules of $A$ . This condition has to be imposed because there are
  examples of QM-abelian surfaces defined over $\mathbb{Q}$ that are not premodular
  (and a fortiori they are not modular, see [DR2]).
   Observe that the action of the quaternion algebra can not be
  defined over $\Q$ because $\Q \subseteq \R$. In fact (we keep from now on
  the premodularity
   condition) the minimal
  field $K$ such that $\End^0_K(A)= B$ is an imaginary quadratic field (cf.
  [DR1]).
  \\ A generalization of the Shimura-Taniyama
conjecture predicts that any abelian variety of $\GL_2$-type over $\Q$
is modular (see [R2]), thus as a particular case modularity of
all premodular QM-surfaces over $\Q$ is expected.\\
 Examples
of this kind of surfaces have been constructed and studied by Brumer
 (see [HHM]), Hashimoto and Murabayashi (see [HM]) and Bending.\\
The fact that $K$ is imaginary implies that
 the field $\Q(\sqrt{d})$ has to be real (cf. [R2], Prop. (7.2),
second proof).
  This
result combined with a description of rationality fields
of endomorphisms of QM-abelian
 surfaces (see [DR1])
 tells us that a principally polarized
  premodular QM-surface over $\Q$ can only exist if the
 quaternion algebra $B$ is ``exceptional", i.e., there exists $m \mid D= \disc(B)$
 such that $B= (  \frac{-D,m}{\Q})$. \\
 %and we have $\End^0_\Q(A)= \Q(\sqrt{m})$. \\
 In this article we will prove the following result:

 \begin{teo}
 \label{teo:main} Any premodular QM-surface over $\Q$ having good
 reduction at $3$ is modular.
 \end{teo}

Remarks: J. Ellenberg has obtained, as a consequence of his proof of
Serre's conjecture over $\F_9$ under local conditions at $3$ and $5$, a
much stronger result asserting the modularity of every abelian surface
over $\Q$ having real multiplication whenever the surface has good
ordinary reduction at $3$ and $5$ (see [E]). Our result is not contained in his
 because it improves these local conditions, but it is much more
 restrictive because of the QM-multiplication hypothesis.\\
% In fact, our result
 % can be seen as the indefinite quaternion algebra analog of the main
  % result in [ES].  \\
Another result of modularity of premodular QM-surfaces was also
obtained by Hasegawa, Hashimoto and Momose (see [HHM])
 with the condition that it exists
an odd prime ramifying in the quaternion algebra $B$ such that the
surface has good reduction at it. Thanks to this result, we can assume
in the proof of our theorem that $B$ is unramified at $3$.\\

\section{The proof}
Let $A$ be a premodular QM-abelian surface over $\Q$
with $\End_\Q^0(A)= \Q(\sqrt{d})$, $d>0$,
$\cO$ the ring of integers of $ \Q(\sqrt{d})$ and $B$  the indefinite quaternion
algebra inside $\End_{\overline{\Q}}^0(A)$.
Let $N$ be the product of the primes of bad reduction
 of $A$ and $D= \disc(B)$. From now on we will assume $3 \nmid N \cdot D$.\\
In the following lines, we shall recollect a few
 facts about Galois representations attached to $A$ (see [R1],
[R2], [HHM], [Die] and
 [DR1] for references).\\
 For any rational prime $\ell$ and $\la \mid
\ell$  a prime in $\cO$ the Galois-action  on the $\ell$-adic Tate
module of $A$ (because $\End^0_\Q(A)=\Q(\sqrt{d})$) gives  a representation
$$ \re: \Gal(\overline{\Q}/\Q) \rightarrow \GL_2(\cO_\lambda)$$
which is odd, irreducible, and unramified outside $\ell \cdot N$.
The field $\Q(\sqrt{d})$ being real, we have: $\det{\re}=
\chi$, the $\ell$-adic cyclotomic character (cf. [R2]).\\
Further restrictions on the compatible family of Galois representations $\{ \re\}$
 are imposed by the fact that $A$ has QM. The situation is like in the
 case of a modular form having a single inner twist (studied by Momose, Ribet,
 Papier): the field of definition of the QM-action is a quadratic
 field $K$, more precisely, an imaginary quadratic field unramified outside
  $N$. Let $H$ be the
 Galois group with fixed field $K$. For every $\ell \nmid N \cdot D$, the image
 of $\re|_H$
 is a normal subgroup of index at most $2$ of the full image
 of $\re$. On the other hand, for every $\ell \nmid N \cdot D$
$$ \re|_H: \Gal(\overline{\Q}/K) \rightarrow \GL_2(\Z_\ell)$$
and this restriction is surjective  for almost every $\ell$ (cf. [Oh] and [DR1]).
If we
call (for every $p \nmid N$)
$a_p$ the trace of the image of $\re(\Fr \; p)$, $\la \nmid p N D$, we know that
$\Q(\{ a_p\}) = \Q(\sqrt{d})$ and  the following
 condition is satisfied:
$$ a_p^\gamma = \varphi(p) a_p \quad \quad (2.1)$$
for almost every $p \nmid N$, where $\gamma$ is the order two element in
 $\Gal(\Q(\sqrt{d})/\Q)$ and $\varphi$ is the quadratic character
corresponding to $K/\Q$. Serre proved that this is equivalent to
 $\Q(\{ a_p^2\}) = \Q \subsetneqq \Q(\{a_p
\})$.\\
From this, it follows (cf. [R1] and [DR1])
 that for $\ell \nmid N \cdot D $, $\la \mid
\ell$, the image of $\re$ is contained in the subgroup of
$\GL_2(\cO_\la)$ generated by $\GL_2(\Z_\ell)$ and the diagonal matrix:
$$  \begin{pmatrix}
    a_p & 0 \\
  0 & a_p^{-1}
\end{pmatrix}$$
for some $p \nmid \ell N$ verifying $\la \nmid a_p$,
 $a_p \neq 0$ and $\varphi(p) = -1$ (if such a prime $p$ exists).
By (2.1), these  two last conditions are equivalent to:
$a_p = u_p \sqrt{d}$, $u_p \in
\Z[\frac{1}{2}]$, $u_p \neq 0$.\\
Let $\rer$ be the reduction modulo $\lambda$ of $\re$ (take it semisimple).
 If for all primes
 $p \nmid \ell N$ with $\varphi(p)=-1$
 we have $\la \mid a_p$, instead of the description above we will make use of the fact
 that (after suitable conjugation)
 the image of $\rer$ will be contained in $\GL_2(\F_\ell)$. This follows from
 the fact that  in this case formula (2.1) implies   that all traces are in $\F_\ell$,
  and the representation
 being odd (and $\det(\rer) = \chi$) it is well known that this implies that
 it has a model over
 $\F_\ell$ (cf. for example [W1], [W2] or [R3]). Observe finally that for a residual
  representation the property
 of being modular depends only on the traces, so it is stable under conjugation.\\

Thanks to the above description, we see that (as in the proof of
 modularity of $\Q$-curves by Ellenberg and Skinner, see [ES]) we
are in a comfortable situation because we have ``for free" that the
residual representation $\bar{\rho}_{\tres}$, if absolutely irreducible,
 is modular, where $\tres$ is
a prime in $\cO$ dividing $3$.
%If  $\bar{\rho}_{\tres}|_H $ is also absolutely
 %irreducible,
This follows from the theorem of
Langlands and Tunnell as in Wiles' original proof (see [W2] and also [ES]),
because
% $\bar{\rho}_{\tres}|_H $ has image in $\GL_2(\F_3)$ and
 even if in general the image of $\bar{\rho}_{\tres} $ is not contained in
 $\GL_2(\F_3)$  its projectivisation
 %(we quotient by its center)
 $\mathbb{P}(\bar{\rho}_{\tres})$ has image in $\PGL_2(\F_3)$.\\
 %If the restriction
 %to $H$ is reducible, then we are in the ``dihedral case" where it is well known
 %that the residual representation is modular.\\
We divide the rest of the proof in two cases:\\

1) $\bar{\rho}_{\tres}$ absolutely irreducible:\\
As we already explained, we know that the residual representation $\bar{\rho}_{\tres}$
is modular. If we also assume that the restriction to $\Q(\sqrt{-3})$
remains absolutely irreducible, then all conditions are satisfied to
apply Diamond's generalization of Taylor-Wiles modularity result (see [D], [W2], [TW])
 and conclude that $\rho_{\tres}$
 is modular, and therefore that $A$ is modular.\\
 Thus, we assume that the restriction of $\bar{\rho}_{\tres}$ to $\Q(\sqrt{-3})$
 is absolutely reducible. This implies that we are in the ``dihedral case",
 namely, that the image of $\bar{\rho}_{\tres}$ is contained in the
 normaliser $\mathcal{N}$ of a Cartan subgroup $C$ of $\GL_2(\bar{\F}_3)$,
 but not contained
 in $C$. We also know that the restriction to $\Q({\sqrt{-3}})$ of our
 representation has its image inside $C$. \\
Thus, the composition of $\bar{\rho}_{\tres}$ with the quotient $\mathcal{N}/C$:
$$\Gal(\overline{\Q}/ \Q) \rightarrow \mathcal{N}
 \twoheadrightarrow \mathcal{N} /C \quad \quad (2.2)$$
gives the quadratic character corresponding to $\Q(\sqrt{-3})$. In
particular, this quadratic character ramifies at $3$.\\
 We know that, in general, for any $\ell \nmid N$, the restriction of
 the residual representation $\rer$ to the inertia subgroup at $\ell$
 has only two possibilities (this is a classical result due to
 Raynaud with some restrictions on $\lambda$,
 which follows in general for $\ell > 2$ from the fact that the $\lambda$-adic
 representation is
 Barsotti-Tate, therefore
  crystalline with Hodge-Tate weights $0$ and $1$, via an application
   of Fontaine-Laffaille theory, cf. [FL] and [B], chapter 9):
   $$ \rer|_{I_\ell} \cong
\begin{pmatrix}
  \chi & * \\
  0 & 1
\end{pmatrix} \; \; \mbox{or} \; \;
\begin{pmatrix}
  \psi_2 & 0 \\
  0 & \psi_2^\ell
\end{pmatrix}
 $$
where $\psi_2$ is a fundamental character of level $2$.\\
If we suppose  that $\bar{\rho}_{\tres}|_{I_3}$ acts through level $2$
fundamental characters, the image of $I_3$ gives a cyclic group of
order $4 > 2$, thus
%In our situation we can exclude the second case: in fact, the standard
(this is a standard
trick of Serre and Ribet, see [R3]) it has to be contained in $C$. But
this implies that the quadratic character defined by composition (2.2) should be
unramified at $3$, contradicting the fact that this character
corresponds to $\Q(\sqrt{-3})$.\\
 %to deal with the dihedral case proves that
%if $\rer$ falls in the dihedral case for a prime $\ell \geq 5$ the
%quadratic character $\phi$ obtained by composing $\rer$ with the quotient
%$N/C$ is unramified at $\ell$ because the image of $I_\ell$
% is a cyclic group of order $\ell \pm 1 \geq 4$ therefore contained in  $C$.
% Now, for $\ell =3$, if we assume that $I_3$ acts through fundamental characters
%  of level $2$, then
%  its image is cyclic of order $4$ and we still can apply this argument, so
%  composing $\bar{\rho}_{\tres}$ with the quotient $N/C$ gives a
%  quadratic character $\phi$ unramified at $3$ in such a case, so the residual
%   representation becomes reducible if we restrict to a quadratic field $K$
%   unramified at $3$. But this
% XXXXXXXX is incompatible with the original assumptions that
%  XXXXXXXXX $\bar{\rho}_{\tres}$ is absolutely irreducible and its restriction
%  XXXXXXXXX to $\Q({\sqrt{-3}})$ is reducible.\\
Thus, we can assume that we are in the first case of Raynaud's result:
 $\bar{\rho}_{\tres}|I_3$ is reducible. Then, we use again the fact that the
 representations $\re$ are
  crystalline with Hodge-Tate weights $0$ and $1$
 for every  $\ell \nmid N$, $\ell > 2$, to conclude, applying results
 of Breuil (see [B], chapter 9), that the
 local representation $\rho_{\tres}|_{G_3}$ is reducible (we can apply this result
 of classification of crystalline representations at $p=3$
 because the highest
 Hodge-Tate weight is $w=1$ and so: $3 > w+1$). Therefore, all
 conditions are satisfied to apply the result of Skinner and Wiles (see [SW2]) on
 nearly ordinary deformations of (residual) modular irreducible
 Galois representations to conclude that $\rho_{\tres}$, thus $A$, is
 modular.\\

2) $\bar{\rho}_{\tres}$ absolutely reducible:\\
In this case, we exclude again the second possibility in
Raynaud' s theorem. This time this is automatic (Kronecker-Weber) and
we have:
%If
%$\bar{\rho}_{\tres}$ is reducible it holds:
$$\bar{\rho}_{\tres}^{s.s.} = \epsilon \oplus \epsilon^{-1} \chi$$
where $\chi$ is the $\mod$ $3$ cyclotomic character and $\epsilon$ a
character ramifying only at primes in $N$.\\
Thus, excluded the case of level $2$ fundamental characters, we apply
again the results of Breuil (cf. [B], chapter 9) on crystalline
 representations to conclude
that $\rho_{\tres}|_{G_3}$ is reducible. Now, an application of the result
of Skinner and Wiles (see [SW1]) on
 % nearly ordinary
 deformations of residually
reducible Galois representations proves that $\rho_{\tres}$, thus $A$, is
modular.

\section{An example}
The following genus two curve (taken from a family of curves constructed
by Brumer) has the property that its Jacobian has real
multiplication by $\Q(\sqrt{5})$
defined over $\Q$:
$$ y^2 = x^6 + 47 x^5 +365 x^4 + 865 x^3 +400 x^2 + 38 x -4 $$
The primes of bad reduction are $2, 5$ and $127$. Let us call $A$ its Jacobian.
  This example was considered in [HHM], were it is proved that $A$
  verifies $\End_{\overline{\Q}}^0(A)= \End_{\Q(\sqrt{-10})}^0(A) =B$
  where $B$ is the indefinite rational
  quaternion algebra of discriminant $10$. Computing a few characteristic
  polynomials we see that $\End_\Q(A)= \Z[\sqrt{5}]$. Applying theorem \ref{teo:main}
  we conclude that $A$ is modular. \\
   Remark: The results of [HHM] can not
  be applied to prove the modularity of $A$ because the primes ramifying in
  $B$ are also primes of bad reduction. The result of [E] is also insufficient
  because $A$ has bad reduction at $5$.

\section{Generalization}
For higher dimensional abelian varieties of $\GL_2$-type having large
enough
endomorphism algebras, we  can still use exactly the same argument to
prove modularity. For example, consider the case of an abelian variety
$A$ of dimension $2^n$ defined over $\Q$
 with:\\
  $\End_\Q^0(A)$ a totally real number field $E$ of degree $2^n \quad$ (I)\\
 such that:\\
  $\End_{\overline{\Q}}^0(A)= M_{2^{n-1}}(B)$, with $B/\Q$ an
 indefinite quaternion algebra (II).\\
  Following [R2],
 we know that $\Gal(E/\Q)$
 is an abelian group of exponent $2$ and that the traces $a_p := \tr(\re (\Fr \; p))$
  again verify: $\Q( \{ a_p^2 \} ) = \Q $ and the
 two dimensional Galois representations $\{ \re \}$
 have, for every $\la \nmid D = \disc(B)$, the following property:
  the image of the projectivization of the residual representation
  $\rer$ is contained in $\PGL_2(\F_\ell)$. Therefore, from the argument
  given in section 2 we conclude:

  \begin{teo} Let $A$ be an abelian variety of dimension $2^n$
  verifying (I) and (II) above. Assume that $A$ has good reduction at
  $3$. Then $A$ is modular.
  \end{teo}

Remark: Again, we can assume in the proof that $ 3 \nmid D = \disc(B)$
 (if $3 \mid D$ modularity follows from [HHM], theorem 2.1).\\

 Final Remark: One of the referees suggested that the arguments and
 results in this article can be easily extended to the case of
 multiplicative reduction at $3$. Instead of applying the results of Breuil,
  all the desired
 statements about restrictions of $\re$ to decomposition groups should follow
 in this case from the
 Mumford-Tate uniformization.

\section{Bibliography}

[B] Breuil, C. ``$p$-adic Hodge theory, deformations and local Langlands",
notes of a course at CRM, Bellaterra, Spain, July 2001\\

[C] Clark, P. ``Two Quaternionic modular philosophies", available at:\\
http://modular.fas.harvard.edu/mcs/archive/Fall2001/notes/10-15-01/ \\

[D] Diamond, F. ``On deformation rings and Hecke rings", Ann. Math. 144
(1996)\\

[Die] Dieulefait, L. ``Newforms, Inner Twists, and the Inverse Galois Problem
 for Projective Linear Groups", J.
  Th. des Nombres de Bordeaux 13 (2001)\\

[DR1] Dieulefait, L. and Rotger, V. ``The arithmetic of QM-abelian surfaces
 through their Galois representations",
preprint \\

[DR2] Dieulefait, L. and Rotger, V. ``Non-modular fake elliptic
curves", preprint \\

[E] Ellenberg, J. ``Serre's conjecture over $\F_9$", available at:\\
front.math.ucdavis.edu/math.NT/0107147 \\

[ES] Ellenberg, J. and Skinner, C. ``On the modularity of $Q$-curves",
 Duke Math. J. 109 (2001)\\

[FL]  Fontaine, J.M. and Laffaille, G. ``Construction de repr{\'e}sentations
$p$-adiques", Ann. Scient. Ec. Norm. Sup., $4^e$ s{\'e}rie, 15 (1982)\\

[HHM] Hasegawa, Y., Hashimoto, K. and Momose, F. ``Modularity conjecture
 for $Q$-curves and
QM-curves", Int. J. of Math. 10 (1999)\\

[HM] Hashimoto, K. and Murabayashi, N. ``Shimura curves as
 intersections of Humbert surfaces and defining
equations of QM-curves of genus two", Tohoku Math. J. (2) 47 (1995)\\

[Oh] Ohta, M. ``On $\ell$-adic representations of Galois groups
obtained from certain two dimensional abelian varieties", J. Fac. Sci.
Univ. Tokyo 21 (1974)\\

[R1] Ribet, K. ``On $\ell$-adic representations attached to modular forms II",
Glasgow Math. J. 27 (1985)\\

[R2] Ribet, K. ``Abelian varieties over $\Q$ and modular forms", KAIST
workshop, Taejon, Korea (1992)\\

[R3] Ribet, K.  ``Images of semistable Galois
representations", Pacific J. of Math.  181 (1997)\\

[SW1] Skinner, C. and Wiles, A. ``Residually reducible representations and modular
forms", IHES Publ. Math. 89 (2000)\\

[SW2] Skinner, C. and Wiles, A. ``Nearly ordinary deformations of irreducible residual
representations", Ann. Fac. Sc. Toulouse Math. (6) 8 (2001)\\

[TW]  Taylor, R. and Wiles, A. ``Ring theoretical properties of certain
 Hecke algebras",
 Ann. of Math. 141 (1995)\\

[W1]   Wiles, A. ``On ordinary $\lambda$-adic representations
associated to modular forms",  Invent. Math. 94 (1988) \\

[W2]   Wiles, A. ``Modular elliptic curves and Fermat's Last Theorem", Ann. of Math.
 141 (1995)

\end{document}